\documentclass [12pt] {report}
\usepackage[dvips]{graphicx}
\usepackage{epsfig, latexsym,graphicx}
\usepackage{amssymb}
\newcommand {\D}[2] {\displaystyle\frac{\partial{#1}}{\partial{#2}}}

\newcommand {\al} {\alpha}

\newcommand {\ga} {\gamma}

\newcommand {\si} {\sigma}

\newcommand {\de} {\delta}
\newcommand {\prtl} {\partial}
\newcommand {\fr} {\displaystyle\frac}

\newcommand {\be} {\begin{equation}}
\newcommand {\ee} {\end{equation}}
\newcommand {\ba} {\begin{array}}
\newcommand {\ea} {\end{array}}
\newcommand {\bp} {\begin{picture}}
\newcommand {\ep} {\end{picture}}
\newcommand {\bc} {\begin{center}}
\newcommand {\ec} {\end{center}}
\newcommand {\bt} {\begin{tabular}}
\newcommand {\et} {\end{tabular}}
\newcommand {\lf} {\left}
\newcommand {\rg} {\right}

\newcommand {\cF} {{\cal F}}

\newcommand {\cR} {{\cal R}}

\newcommand {\ses} {\medskip}

\newcommand {\e} {\mathop{\rm e}\nolimits}

\newcommand {\bibit} {\bibitem}
\newcommand {\nin} {\noindent}

\usepackage{psfrag}

\usepackage{epsfig, latexsym,graphicx}
 \usepackage{amsmath}

\newcommand {\cD} {{\cal D}}

\def\2#1#2#3{{#1}_{#2}\hspace{0pt}^{#3}}
\def\3#1#2#3#4{{#1}_{#2}\hspace{0pt}^{#3}\hspace{0pt}_{#4}}
\newcounter{sctn}
\def\sec#1.#2\par{\setcounter{sctn}{#1}\setcounter{equation}{0}
                  \noindent{\bf\boldmath#1.#2}\bigskip\par}

\begin {document}

\begin {titlepage}

\vspace{0.1in}

\begin{center}

{\Large \bf  Finsleroid-regular   space.    Landsberg-to-Berwald implication}

\end{center}

\vspace{0.3in}

\begin{center}

\vspace{.15in} {\large G.S. Asanov\\} \vspace{.25in}
{\it Division of Theoretical Physics, Moscow State University\\
119992 Moscow, Russia\\
{\rm (}e-mail: asanov@newmail.ru{\rm )}} \vspace{.05in}

\end{center}

\begin{abstract}

\ses

By performing  required   evaluations, we show that
in the Finsleroid-regular   space the Landsberg-space condition just degenerates to the Berwald-space
condition (at any dimension number $N\ge2$).
Simple and clear expository representations are obtained.
Due comparisons   with the Finsleroid-Finsler space  are indicated.

\ses


\end{abstract}

\end{titlepage}

\vskip 1cm

\ses

\ses

\setcounter{sctn}{1} \setcounter{equation}{0}

\nin
  {\bf 1. Description  of new conclusions}

\ses

\ses

The Finsler geometry theory can mentally be divided into two great parts:
the development of the general theory to search for tensorial and geometrical implications
 of a general Finsler metric function ``$F$''
and the investigation of possible results of specifying the $F$ in an attractive particular way
[1,2].
Obviously, it is the latter way that one is to follow when hoping to develop genius and handy applications.

Below, we   deal with the Finsler space notion specified by the condition
that the basic  Finslerian metric function, to be denoted by  $K(x,y)$,
 is constructed functionally from the set
$\{g(x), b_i(x), a_{ij}(x),y\}$,
where
$g(x)$ is a scalar,
$b_i(x)$ is an involved vector field,
 and $a_{ij}(x)$ is a  (positive-definite) Riemannian metric tensor; $y$ stands for the tangent vectors supported by a
 point $x$ of the underlined manifold.
Denoting  by
\be
c=||b|| \equiv || b||_{\rm{Riemannian}}
\ee
the respective Riemannian norm value of the input 1-form $b=b_i(x)y^i$
and assuming the range
\be
 0 < c < 1,
 \ee
we construct the particular function $K(x,y)$  which occurs being globally regular.
The  entailed positive-definite Finsler space  will be denoted by ${\mathbf\cF\cR^{PD}_{g;c} } $.

The extrapolation
\be
\cF\cR^{PD}_{g;c\to 1}    =   \cF\cF^{PD}_g
\ee
takes place, where $\cF\cF^{PD}_g$ is the Finsleroid-Finsler space which was constructed and developed
in [3-6] under the assumption $||b||=1$.

This  scalar    $c(x)$ proves to play the role of the {\it regularization factor.}
Indeed, in the spaces ${\mathbf\cF\cR^{PD}_{g;c} } $ and $ \cF\cF^{PD}_g$
% (as well as in the space   $ \cF\cF^{SR}_g$)
the metric function $K$   is  constructed such that $K$  involves
 the square-root variable
$q(x,y)=\sqrt{S^2-b^2}$ (see (A.4)).
Differentiating various  tensors of  the space ${\mathbf\cF\cR^{PD}_{g;c} } $ as well as the space
 $ \cF\cF^{PD}_g$ gives rise,
therefore,
 to appearance of degrees of the fraction $1/q$.

If $||b||=1$, we have  $q=0$  when  $y=\pm b$.

As far as $ 0<c<1$, we have $q\ne 0$ whenever $y\ne 0$
(because of the inequality (A.5)), so that the fraction $1/q$
 does not  produce any singularities on $TM\setminus 0$.

We use

\ses

\ses

{REGULARITY   DEFINITION.}   The Finsler space  $\cF\cR^{PD}_{g;c} $
under consideration is {\it regular} in the following sense: {\it globally}
over all the slit tangent bundle $TM\setminus 0$,
the Finsler metric function $K(x,y)$ of the space is smooth of the class $C^{\infty}$
regarding both the arguments $x$ and $y$,
and also the entailed Finsler metric tensor $g_{ij}(x,y)$ is positive-definite: $\det( g_{ij}) >0$.

\ses

\ses

{%\pgbrk%}

The $\cF\cF^{PD}_g $--space  is  smooth  of the class  $C^2$, and not of the class $C^3$,
  on all of the  {\it slit tangent bundle }
   $TM\setminus 0$.
The $\cF\cF^{PD}_g $--space  is  smooth  of the class  $C^{\infty}$
on all of the {\it $b$-slit  tangent bundle }
\be
{\cal T}_bM~:= TM \setminus0\setminus b\setminus -b
\ee
(obtained
by deleting out in $TM\setminus 0$
all the directions which point along, or oppose, the directions  given rise to by the  1-form $b$).

The Finsleroid-Finsler space $ \cF\cF^{PD}_g$ developed
 involves an attractive realization of the
Landsberg condition over the  $b$-slit  tangent bundle (see [3-6]).
The realization cannot be extrapolated to the
$b$-section of the  tangent bundle $TM$, because on the section  the smoothness of the space
 $\cF\cF^{PD}_g $
degenerates to but the $C^2$-level.

It is impossible to lift realization of the Landsberg condition from the space
$\cF\cF^{PD}_g $
to the space
${\mathbf\cF\cR^{PD}_{g;c} } $. Indeed, the following theorem is valid.

\ses

{\bf Landsberg-to-Berwald Theorem.}
In  the Finsleroid-regular space $\cF\cR^{PD}_{g;c}$ the Landsberg-space condition entails the Berwald-space condition.

\ses

To verify this theorem,  it is sufficient to pay a due attention to the factor $(1-c^2)$
which enters the right-hand part of the formula (A.32) which precedes  the explicit (and simple)
expression (A.33) for the contraction
$A^iA_i$ (see more detail in Note placed in the end of Appendix A).

\ses

In Appendix A,
the explicit form of the Finsleroid-regular metric function $K$ is presented,
the space ${\mathbf\cF\cR^{PD}_{g;c} } $ is rigorously defined,
and entailed representations of various key   tensors are given.
The knowledge of the associated spray coefficients $G^i$,
as given by the explicit representation (A.37) derived,
can open up various convenient possibilities to evaluate and study associated geodesic equations, connection coefficients,
as well as curvature tensors.
Amazingly, the coefficients (A.37) provide us readily with the  Berwald space of the regular type:
the Berwald case would imply $g=const$ and $\nabla_ib_j=0$ (see (A.35) and (A.36)), and also $c=const$.
In the two-dimensional case, however, the space degenerates to the locally Minkowskian space (whenever $g\ne0$).
In the dimensions $N\ge 3$ the obtainable ${\mathbf\cF\cR^{PD}_{g;c} } $-Berwald spaces can be
``neither Riemannian  nor locally Minkowskian."
 The Berwald case of the Finsler space is attractive because of its simplicity.
\ses

{%\pgbrk}

In Appendix B,
we investigate the derivatives of the spray coefficients in the particular case (B.1).
Calculations involved are simple, showing the validity of the following theorem.
\ses

\ses

{\bf Particular Theorem}.
Given a scalar
$k=k(x)$.
If the conditions   $   g=const $ and $  \nabla_ib_j=kr_{ij}$   are fulfilled,
then
\be
\dot A_{knj}=
(1-c^2)k
\bigl(m_1A_{knj}
+m_2A_kA_nA_j
\bigr).
\ee

 \ses

 \ses

In (1.5),  $m_1$ and $m_2$ are two scalars, which form is indicated explicitly in (B.20) of Appendix B.
Remarkably,
$m_1$, as well as $m_2$, doesn't vanish identically when $g\ne0$.
Therefore, if $g\ne 0$ and $0<c<1$,
the right-hand part in (1.5) can vanish identically only in the case $k=0$ which is the Berwald case.

We shall use the notation
\be
\cD=\fr1K   y^jD_j
\ee
 with $D_j$ standing for the
 $h$-covariant (horizontal) Finslerian derivative (see [2]);
the tensor $\dot A_{knj}$ is identical to that used in [2],
namely,
\be
\dot A_{knj}~:=\cD   A_{knj}.
\ee

The  condition $\nabla_ib_j=kr_{ij}$, when taken in conjunction with $g=const$, realizes
 the  Landsberg space, that is, $\dot A_{knj}=0$,
  in the Finsleroid-Finsler space
$\cF\cF^{PD}_g$
(see [3-6]). Lifting the condition to the Finsleroid-regular space
${\mathbf\cF\cR^{PD}_{g;c} } $ results in the representation
(1.5) which (beautifully?)  extends the Landsberg condition $\dot A_{knj}=0$.

\ses

The occurrence of the regularizing  factor $(1-c^2)$  in the right-hand parts of the formulae  (A.32)
and
(1.5)
 presents an astonishingly simple and  explicit illustration to  the above
 Landsberg-to-Berwald Theorem.

\ses

{%\pgbrk}

\ses\ses

\setcounter{equation}{0}

{
\nin \bf Appendix A:   Involved  ${\mathbf\cF\cR^{PD}_{g;c} } $-notions
 }

\ses\ses

Let $M$ be an $N$-dimensional
$C^{\infty}$
differentiable  manifold, $ T_xM$ denote the tangent space to $M$ at a point $x\in M$,
and $y\in T_xM\backslash 0$  mean tangent vectors.
Suppose we are given on $M$ a Riemannian metric ${\cal S}=S(x,y)$.
 Denote by
$\cR_N=(M,{\cal S})$
the obtained $N$-dimensional Riemannian space.
Let us also assume that the manifold $M$ admits a non--vanishing 1-form $b= b(x,y)$,
denote by
\be
c=||b|| \equiv || b||_{\rm{Riemannian}}
\ee
the respective Riemannian norm value.
Assuming
\be
 0 < c < 1,
 \ee
we get
\be
S^2-b^2>0
\ee
and may conveniently use the variable
\be
q:=\sqrt{S^2-b^2}.
\ee
Obviously, the inequality
\be
q^2 \ge \fr{1-c^2}{c^2}\,b^2
\ee
is valid.

With respect to  natural local coordinates in the space
$\cR_N$
we have the local representations
\be
\sqrt{a^{ij}(x)b_i(x)b_j(x)} =    c(x)
\ee
and
\be
 b=b_i(x)y^i,  \qquad S= \sqrt{a_{ij}(x)y^iy^j}.
\ee
The reciprocity  $a^{in}a_{nj}=\de^i{}_j$ is assumed, where $\de^i{}_j$ stands for the Kronecker symbol.
The covariant index of the vector $b_i$  will be raised by means of the Riemannian rule
$ b^i=a^{ij}b_j,$ which inverse reads $ b_i=a_{ij}b^j.$
We  also  introduce the tensor
\be
r_{ij}(x)~:=a_{ij}(x)-b_i(x)b_j(x)
\ee
  to have the representation
\be
  q=\sqrt{r_{ij}(x)y^iy^j}.
  \ee

{%\pgbrk}

We choose  the Finsler space notion specified by the condition
that the Finslerian metric function
$K(x,y)$ be of the functional dependence
\be
K(x,y) =\Phi \Bigl(g(x), b_i(x), a_{ij}(x),y\Bigr),
\ee
where $g(x)$ is a scalar
(on the background manifold $M$),
 subjected to ranging
\be
-2<g(x)<2,
\ee
and apply  the convenient notation
\be
h(x)=\sqrt{1-\fr14(g(x))^2}, \qquad
G(x)=\fr{g(x)}{h(x)}.
\ee
We introduce
the {\it  characteristic  quadratic form}
\be
B(x,y) :=
b^2+gqb+q^2
\equiv
\fr 12
\Bigl[(b+g_+q)^2+(b+g_-q)^2\Bigr],
\ee
where $ g_+=(1/2) g+h$ and $ g_-=(1/2) g-h$,
The  discriminant $D_{\{B\}}$ of the quadratic form $B$
is negative:
\be
D_{\{B\}}   =    -4h^2<0.
\ee
Therefore, the {\it quadratic form $B$ is  positively definite.}
In the limit $g\to 0$,
the definition (A.13) degenerates to the
 quadratic form  of the input Riemannian metric tensor:
$   B|_{_{g=0}}=b^2+q^2 \equiv S^2.   $
Also,
$  \eta B|_{_{y^i=b^i}}=c^2,  $
where
$ \eta=1/(1+gc\sqrt{1-c^2}).  $
It can readily be verified that on the definition range (A.11) of the $g$  we have
$   \eta>0.   $

Under these conditions,
we set forth  the following definition.

\ses

\ses

 {\large  Definition}. The scalar function $K(x,y)$ given by the formulas
\be
K(x,y)=
\sqrt{B(x,y)}\,J(x,y)
\ee
and
\be
J(x,y)=\e^{-\frac12G(x)f(x,y)},
\ee
where
\be
f=
-\arctan
 \fr G2
+\arctan\fr{L}{hb},
\qquad  {\rm if}  \quad b\ge 0,
\ee
and
\be
f= \pi-\arctan
\fr G2
+\arctan\fr{L}{hb},
\qquad  {\rm if}
 \quad b\le 0,
\ee
 with
 \be
 L =q+\fr g2b,
\ee
\ses\\
is called
the {\it  Finsleroid-regular   metric function}.

\ses

\ses

The function $L$  obeys   the identity
\be
L^2+h^2b^2=B.
\ee

\ses

\ses

 {\large  Definition}.  The arisen  space
\be
\cF\cR^{PD}_{g;c} :=\{\cR_{N};\,b_i(x);\,g(x);\,K(x,y)\}
\ee
is called the
 {\it Finsleroid-regular  space}.

\ses

\ses

 {\large  Definition}. The space $\cR_N$ entering the above definition is called the {\it associated Riemannian space}.

\ses

\ses

The associated Riemannian metric tensor $a_{ij}$ has the meaning
\be
 a_{ij}=g_{ij}\bigl|_{g=0}\bigr. .
\ee

\ses

\ses

 {\large  Definition}.  Within  any tangent space $T_xM$, the Finsleroid-regular metric function $K(x,y)$
 produces the {\it regular Finsleroid}
 \be
 \cF\cR^{PD}_{g;c\, \{x\}}:=\{y\in   \cF\cR^{PD}_{g;c\, \{x\}}: y\in T_xM , K(x,y)\le 1\}.
  \ee

\ses

 \ses

 {\large  Definition}. The {\it regular Finsleroid Indicatrix}
 $ I\cR^{PD}_{g;c\, \{x\}} \subset T_xM$ is the boundary of the regular Finsleroid, that is,
 \be
 I\cR^{PD}_{g;c\, \{x\}} :=\{y\in I\cR^{PD}_{g;c\, \{x\}} : y\in T_xM, K(x,y)=1\}.
  \ee

\ses

 \ses

 {\large  Definition}. The scalar $g(x)$ is called
the {\it Finsleroid charge}.
The 1-form $b=b_i(x)y^i$ is called the  {\it Finsleroid--axis}  1-{\it form}.

\ses

\ses

We shall meet the function
\be
\nu:=q+(1-c^2)gb
\ee
for which
\be
\nu >0 \quad \rm{when} \quad |g|<2.
\ee
Indeed, if $gb>0$, then  the right-hand part of (A.25) is positive.
 When $gb<0$, we may note that
 at any fixed $c$ and $b$ the minimal value of $q$ equals
$\sqrt{1-c^2}|b|/c$ (see (A.5)), arriving again at (A.26).

Under these conditions,
 we can   explicitly extract from the function $K$ the  distinguished Finslerian tensors,
 and first of all
the covariant tangent vector $\hat y=\{y_i\}$ from $y_i :=(1/2)\partial {K^2}/ \partial{y^i}$,
obtaining
\be
y_i=(u_i+gqb_i) \fr{K^2}B,
\ee
where $u_i=a_{ij}y^j$.
After that, we can find
the  Finslerian metric tensor $\{g_{ij}\}$
together with the contravariant tensor $\{g^{ij}\}$ defined by the reciprocity conditions
$g_{ij}g^{jk}=\de^k_i$, and the  angular metric tensor
$\{h_{ij}\}$, by making  use of the following conventional  Finslerian  rules in succession:
$$
g_{ij} :
=
\fr12\,
\fr{\prtl^2K^2}{\prtl y^i\prtl y^j}
=\fr{\prtl y_i}{\prtl y^j}, \qquad
h_{ij} := g_{ij}-y_iy_j\fr1{K^2},
$$
thereafter   the Cartan tensor
\be
 A_{ijk}~ := \fr K2\D{g_{ij}}{y^k}
 \ee
and the contraction
\be
 A_i~:=g^{jk}A_{ijk} =  K\D{\ln\bigl(\sqrt{\det(g_{mn})}\bigr)}{y^i}
\ee
can readily be evaluated.

{%\pgbrk}

It can straightforwardly be verified that
\be
\det(g_{ij})=\fr {\nu}q\biggl(\fr{K^2}B\biggr)^N\det(a_{ij})>0
\ee
with  the function $\nu$ given by (A.25) [7],
and
\be
A_i=\fr {Kg}{2qB}
\fr1X
(q^2b_i- bv_i),
\ee
where
 the function $X$ is given by
\be
\fr1X=    N+    (1-c^2)     \fr{B}{q\nu}.
\ee
Contracting yields the formula
\ses
\be
A^iA_i=     \fr{g^2}4
\fr1{X^2}\lf(N+1-\fr1X\rg)
\ee
and
 evaluating  the Cartan tensor results in the lucid representation
\be
A_{ijk}= X
 \Biggl[
A_ih_{jk}  +A_jh_{ik}  +A_kh_{ij}
-\lf(N+1-\fr1X\rg)
\fr1{A_hA^h}A_iA_jA_k
\Biggr].
\ee

\ses

\ses

We use
 the Riemannian covariant derivative
\be
\nabla_ib_j~:=\partial_ib_j-b_ka^k{}_{ij},
\ee
where
\be
a^k{}_{ij}~:=\fr12a^{kn}(\prtl_ja_{ni}+\prtl_ia_{nj}-\prtl_na_{ji})
\ee
are the
Christoffel symbols given rise to by the associated Riemannian metric ${\cal S}$.

%%%%%%   figure           Ffigure

{%\pgbrk}

Attentive direct calculations
of
the induced spray coefficients
$ G^i =\ga^i{}_{nm}y^ny^m$,
where
$\ga^i{}_{nm}$ denote the associated Finslerian Christoffel symbols,
can be used to arrive at the following result.

\ses
\ses

 {\large Theorem 1.}  {\it In the Finsleroid-regular  space
 ${\mathbf\cF\cR^{PD}_{g;c} } $
the spray coefficients $G^i$
can  explicitly be written in the form}
\be
G^i=
\fr  g{\nu}
\Bigl(
y^jy^h\nabla_jb_h
+gqb^jf_j\Bigr)
v^i
-gqf^i
 +E^i
 +a^i{}_{nm}y^ny^m.
\ee

\ses

\ses

We use  the  notation
\be
v^i=y^i-bb^i
\ee
and
\be
f_j=f_{jn}y^n,\quad
f^i=f^i{}_ny^n,\quad
f^i{}_n=a^{ik}f_{kn}, \quad
f_{mn}=
\nabla_mb_n-\nabla_nb_m
\equiv \D{ b_n}{x^m}-\D {b_m}{x^n},
\ee
where
  $\nabla$ means the covariant derivative in terms of the associated Riemannian space
${R}_N=(M,{S})$  (see (A.35));
  $ a^i{}_{nm} $ stands for
the  Riemannian  Christoffel symbols (A.36) constructed from the input Riemannian metric tensor $a_{ij}(x)$;
the coefficients $E^i$ involving the gradients
 $g_h=\partial g/\partial x^h$
 of the Finsleroid charge
 can be taken as
\be
E^i = {\bar M} (yg)y^i
+  K \fr{2 q^2 }{gB}(yg) X A^i
-\fr12 {\bar M}K^2g_hg^{ih},
\ee
where $(yg)=g_hy^h$,
$X$ is the function given in (A.32), and
the function ${\bar M}$  is defined by the equality $\partial K/\partial g= (1/2){\bar M}K$.
It can readily be verified that
\be
b_hE^h=b{\bar M}(yg)
+
\Biggl[  \fr{q^2}{B}(yg)
 -  \fr12  {\bar M}
 \Bigl(
 q(bg)
+ g[b(bg)-(yg)]
 \Bigr)
 \Biggr]
  \fr {1}{\nu}  (c^2S^2-b^2)
  -\fr12 {\bar M}(bg)b^2,
\ee
where $(bg)=b^hg_h$.

\ses

Contracting (A.37) by $b_i$ yields
$$
b_i
\Bigl(G^i-a^i{}_{mn}y^my^n\Bigr)
=
\fr  g{\nu}
\Bigl(
(ys)
+gq \si\Bigr)
(1-c^2)b
-gq\si
 +b_iE^i,
$$
\ses
or
\be
b_k
\Bigl(G^k-a^k{}_{mn}y^my^n\Bigr)
=
\fr  g{\nu}
(ys)
(1-c^2)b
-\fr  {gq^2}
{\nu}\si
 +b_kE^k,
\ee
where    $(ys)=y^hy^k \nabla_hb_k$   and
$$
\si=b_kf^k=b^ky^nf_{kn}.
$$

Let us apply  the operator
$  \cD  $
(defined in (1.6)) to the input 1-form $b$:
$$
\cD b=y^k\partial_k b-G^kb_k
=y^ky^m\nabla_kb_m   -(G^k-a^k{}_{mn}y^my^n)b_k,
$$
that is,
\be
\cD b=(ys)-(G^k-a^k{}_{mn}y^my^n)b_k.
\ee
\ses
Taking into account   (A.42) and denoting   $\dot b=\cD b$,
we obtain simply
\be
\dot b=
\fr  q{\nu}
(ys)
+\fr  {gq^2}
{\nu}\si
-b_kE^k,
\ee

\ses

{%\pgbrk}

With the help of the coefficients
\be
G^k{}_n=\D{G^k}{y^n}
\ee
we can consider also the contracted derivative
\be
\cD b_n=
y^h\partial_hb_n
-
\fr12 G^k{}_n b_k
=
y^h\nabla_hb_n   -   \lf(\fr12 G^k{}_n-a^k{}_{nh}y^h\rg)  b_k,
\ee
or
\be
\cD b_n=
y^h\nabla_hb_n
-\fr12
\D{\Biggl(b_k\Bigl(G^k-a^k{}_{mn}y^my^n\Bigr)\Biggr)}
{y^n}.
\ee
This, by virtue of (A.43), can be written as
\be
\cD b_n=
\fr12
y^h(\nabla_hb_n-\nabla_nb_h)
+\fr12
\D{\cD b}
{y^n}.
\ee

In evaluations, it is convenient to use the derivative value
\be
\nu_k=\fr{v_k}q+(1-c^2)gb_k \equiv \D{\nu}{y^k}.
\ee

If the input 1-form $b=y^ib_i$ is exact:
\be
f_{mn}=0,
\ee
then the above representation (A.37) reduces to read merely
\be
G^i=
\fr  g{\nu}
(ys)
v^i
+E^i
+a^i{}_{nm}y^ny^m,
\ee
entailing
\be
G^i{}_k=
gU_kv^i
+g\fr1{\nu} (ys)
r^i{}_k
+E^i{}_k
+2a^i{}_{km}y^m
\equiv
\D{G^i}{y^k}
\ee
with
\be
U_k=-\fr1{\nu^2}\nu_k
(ys)
+\fr2{\nu} s_k,
\ee
\ses
where
$(ys)=y^hs_h$ and $s_k=y^h\nabla_hb_k$;
~
$E^i{}_k=\partial E^i/ \partial y^k$
and  $r^i{}_k=a^{ih}r_{hk}$.
Differentiating (A.53) yields the coefficients
$U_{km}=\partial U_{k} / \partial y^m$ given by the formula
\be
U_{km}=2\fr1{\nu^3}\nu_k\nu_m(ys)
-\fr2{\nu^2}
(\nu_k s_m+\nu_m s_k)
-\fr1{\nu^2q}\eta_{km} (ys)
+\fr2{\nu}\nabla_mb_k,
\ee
where
\be
\eta_{km}=r_{km}-\fr1{q^2}v_kv_m, \qquad v_k=u_k-bb_k \equiv a_{kj}v^j, \qquad u_k=a_{kj}y^j.
\ee

{%\pgbrk}

From (A.52) we find that
 the coefficients
$  G^i{}_{km}=  \partial G^i{}_{k}/ \partial y^m$
are given by the representation
\be
G^i{}_{km}=
gU_kr^i{}_m
+
gU_mr^i{}_k
+gU_{km}v^i
+E^i{}_{km}
+2a^i{}_{km},
\ee
where
$  E^i{}_{km}=  \partial E^i{}_{k}/ \partial y^m$.
By contracting
we obtain
\be
b_iG^i{}_{km}=(1-c^2)\Bigl(gU_kb_m+gU_mb_k+gU_{km}b\Bigr)
+b_iE^i{}_{km}
+2b_ia^i{}_{km}
\ee
and
\be
u_iG^i{}_{km}=\Bigl(gU_kv_m+gU_mv_k+gU_{km}q^2\Bigr)
+u_iE^i{}_{km}
+2u_ia^i{}_{km}.
\ee

{%\pgbrk}

\ses

The following theorem is valid.
\ses
\ses

 {\large  Theorem 2.}
   {\it If the Finsleroid-regular  space
 ${\mathbf\cF\cR^{PD}_{g;c} } $  is the  Berwald space and is not the locally-Minkowskian space,
 then
$g=const$.
}

\ses

\ses

To verify the theorem, it is worth noting  that in the Berwald case the covariant derivative
$\cD_n A_{ijk}$ of the Cartan tensor
 vanishes identically (see [1,2]). The implication
 $\cD_n A_{ijk}=0 ~ \Longrightarrow \cD_n (A^iA_i)=0$ is obviously valid in any Finsler space.
Therefore,  in the Berwald case of dimension $N\ge3$  the representations  (A.33) and (A.34) just entail $g=const$.
In the two-dimensional case the representation (A.34) reduces to
\be
A_{ijk}=
I  \al_i \al_j \al_k    \quad {\text at} ~~  N=2
\ee
with   the normalized vector   $\al_i=A_i/\sqrt{A^hA_h}$,  and with the main scalar $I$ given by
\be
I=   \sqrt{A^hA_h}.
\ee
It is well-known [1,2]
that the two-dimensional Finsler space is the non-Minkowskian Berwald space if and only if the main scalar is
independent of the argument $y$. However, from (A.32) and (A.33) it just follows that
the identical vanishing of the derivative   $\partial (A^hA_h) / \partial y^n $ entails $g=0$,
 that is, the Riemannian space. Thus, the theorem  is valid  in any dimension $N\ge2$.

\ses

Also, the following theorem is valid.

\ses
\ses

 {\large  Theorem 3.}
   {\it The Finsleroid-regular  space
 ${\mathbf\cF\cR^{PD}_{g;c} } $
is the Berwald space if and only the conditions
\be
g=const ~ ~ { \rm and}  ~  ~
\nabla_mb_n=0
\ee
hold.}

\ses

\ses

The vanishing
$\nabla_mb_n=0$  can well be interpreted geometrically by phrasing that the 1-form $b$
is {\it parallel} (in the sense of the associated Riemannian space).

When  the coefficients
 $
                                                    %G^i{}_{nm}=
(1/2)\partial^2G^i/\partial y^m \partial y^n$
 are independent of the variable $y$,
 one says that the Finsler space is the Berwald space (see [1,2]).
The sufficiency of the conditions (A.61) is obvious from the representation (A.37),
reducing the coefficients $(1/2)\partial^2G^i/\partial y^m \partial y^n$
 to
the Riemannian Christoffel symbols:
 \be
 G^i=  a^i{}_{nm}y^ny^m ~ ~ \text {in the Berwald case}.
\ee
To verify the necessity, we apply Theorem 2 to conclude
 $g=const$,
which in turn yields $E^i=0$ in the representation (A.37) of $G^i$,
after which
it is easy to see that the Berwald space arises  if only  $\nabla_mb_n=0$,
as far as the value of $g$ is kept differing from zero (the choice $g=0$ would reduce the Finsler
 spaces  under consideration to  Riemannian spaces).

The above theorem yields an attractive example of the regular Berwald space.

{%\pgbrk}

\ses

Comparing the conditions (A.61) with the representation (A.44) of $\dot b$ yields the following elegant theorem.

\ses
\ses

 {\bf Notable Theorem.}
   {\it In the dimensions $N \ge3$,
   the Finsleroid-regular  space
 ${\mathbf\cF\cR^{PD}_{g;c} } $
is the non-Riemannian Berwald space if and only if
}
\be
\dot b=0  ~~ {\text and} ~~ g=const\ne0.
\ee

\ses
\ses

From (A.44) it is obvious that  (A.61) entails  $\dot b=0$.
The opposed implication that  the conditions $\dot b=0$
 and $ g=const\ne0$ entail
 $\nabla_mb_n=0$
  can straightforwardly and readily be arrived at on the basis of the representation (A.44).

\ses

\ses

NOTE.
The Landsberg-space condition means the requirement that the (identical) vanishing
$\cD A_{ijk}=0$ hold.
Obviously,
$\cD A_{ijk}=0 ~ \Longrightarrow \cD (A^iA_i)=0$.
In dimensions $N\ge3$, the Cartan tensor representation (A.34)  communicates us that
$\cD A_{ijk}=0$ would entail $\cD X=0$. When
$\cD (A^iA_i)=0$ and $\cD X=0$  are applied to the representation (A.33) of $A^iA_i$, we conclude
$\cD g=0$. Since $g$ is a function of $x$, the last vanishing entails  $g=const$.
In the two-dimensional case,
 the Cartan tensor representation (A.34) reduces to the representation (A.59)-(A.60) which doesn't involve the function
 $X$, which circumstance does not make possible  to conclude $\cD X=0$ in a direct way. However,
the basic definition of $K$, and hence the derivative of $K$ with respect to $g$,
  involves the function $\arctan(L/hb)$ (see (A.17)).
  Obviously, this trigonometric  function cannot be cancelled by polynomials constructed from the
variable set $\{b,q\}$. At the same time, simple evaluations show that the  function enters the
quantity $\cD X$  through the term $(1-c^2)f(x,y) (\arctan(L/hb))\cD g$ with some function $f$
which doesn't vanish identically. The identical vanishing of   $\cD g$ [which is equal to
 $ (yg)=y^j\partial g/\partial x^j$]   means $g=const$.
Therefore,
the Landsberg-space condition entails $g=const$  and $\cD X=0$ in any case of the dimension $N\ge 2$.
 Attentive consideration can be applied to conclude after direct  calculations
that whenever $g=const$ the identical vanishing $\cD X=0$ is possible when either $g=0$ (which is the Riemannian case)
or
 $\nabla_mb_n=0$ (which is the Berwald case). In this way  the
 Landsberg-to-Berwald Theorem set forth in Section 1 is getting  valid.

\ses

More detail of calculations involved in the space
$\cF\cR^{PD}_{g;c} $
can be found in [7].

{%\pgbrk}

\ses\ses

\setcounter{equation}{0}

{  \nin \bf Appendix B:   Evaluation of $\dot A_{knj}$      in particular case    }

\ses

\ses

Below we evaluate the particular case
\be
g=const \qquad  {\rm{and}}  \qquad \nabla_ib_j=kr_{ij}, \quad k=k(x),
\ee
obtaining
the following representations from (A.53)-(A.55):
$$
U_k=-k\fr{q^2}{\nu^2}\nu_k
+k\fr2{\nu}v_k,
$$
\ses
$$
G^i{}_k=
gU_kv^i
+gk\fr{q^2}{\nu}
r^i{}_k
+2a^i{}_{km}y^m,
$$
and
\be
U_{km}=2k\fr{q^2}{\nu^3}\nu_k \nu_m
-k\fr2{\nu^2}
(\nu_k v_m+\nu_mv_k)
-k\fr q{\nu^2}\eta_{km}
+k\fr2{\nu}r_{mk}.
\ee

Let us find the tensor
$U_{kmj}=\partial U_{km} / \partial y^j$:
$$
U_{kmj}=-6k\fr{q^2}{\nu^4}    \nu_k   \nu_m   \nu_j
+4k\fr1{\nu^3}   (\nu_k\nu_mv_j   +   \nu_j\nu_mv_k   +   \nu_k\nu_jv_m   )
$$
\ses
$$
+2k\fr q{\nu^3}
(
\nu_m\eta_{kj}
+
\nu_k\eta_{mj}
+
\nu_j\eta_{km}
)
-
2k\fr 1{\nu^2}
(
\nu_mr_{kj}
+
\nu_kr_{mj}
+
\nu_jr_{km}
)
-k\fr 1{q\nu^2}(\eta_{km}v_j+\eta_{jm}v_k+\eta_{kj}v_m),
$$
\ses
\ses
or
$$
U_{kmj}=-6k\fr{q^2}{\nu^4}    \nu_k   \nu_m   \nu_j
+4k\fr1{\nu^3}   (\nu_k\nu_mv_j   +   \nu_j\nu_mv_k   +   \nu_k\nu_jv_m   )
+2k\fr{ q-\nu}{\nu^3}
(
\nu_m\eta_{kj}
+
\nu_k\eta_{mj}
+
\nu_j\eta_{km}
)
$$
\ses
\be
-
2k\fr 1{q^2\nu^2}
(
\nu_mv_kv_j   +   \nu_k  v_mv_j  +   \nu_j  v_kv_m
)
-k\fr 1{q\nu^2}(\eta_{km}v_j+\eta_{jm}v_k+\eta_{kj}v_m).
\ee

After that, we can evaluate the coefficients
$G^i{}_{kmj}=\partial G^i{}_{km}/ \partial y^j$
with the help of (A.56), which yields
 the representation
\be
G^i{}_{kmj}=
gU_{kj}r^i{}_m
+
gU_{mj}r^i{}_k
+gU_{km}r^i{}_j
+gU_{kmj}v^i,
\ee
from which we can find the contraction
\be
u_iG^i{}_{kmj}=
gU_{kj}v_m
+
gU_{mj}v_k
+gU_{km}v_j
+gU_{kmj}q^2.
\ee

\ses

{%\pgbrk}

It is convenient to use the vector
\be
e_k=\fr b{q^2}v_k-b_k,
\ee
having
\be
\nu_k= \fr 1b(q e_k + \nu b_k)
\ee
on the basis of (A.49).

We obtain readily
\be
U_{k} =   kq^2\fr{2\nu-q}{\nu^2 b}e_k    +k\fr{q^2}{\nu b}b_k,
\ee
\ses
together with
\be
U_{km}=   k\fr{2\nu-q}{\nu^2}\eta_{km}   +2k(1-c^2)^2\fr{g^2q^2}{\nu^3}e_ke_m
\ee
\ses
which entails
$$
U_{kmj}=
-\fr1{q^2}(v_kU_{mj}+v_mU_{kj}+v_jU_{km})
$$
\ses
\be
+2k(1-c^2)^2\fr{g^2b}{\nu^3}
(\eta_{kj}e_m+\eta_{mj}e_k+\eta_{mk}e_j)
+6k(1-c^2)^3\fr{g^3q^2}{\nu^4}
e_ke_me_j
\ee
\ses
and
$$
G^i{}_{kmj}=
gU_{kj}\eta^i{}_m
+
gU_{mj}\eta^i{}_k
+gU_{km}\eta^i{}_j
$$
\ses
\be
+g
\Biggl[
2k(1-c^2)^2\fr{g^2b}{\nu^3}
(\eta_{kj}e_m+\eta_{mj}e_k+\eta_{mk}e_j)
+6k(1-c^2)^3\fr{g^3q^2}{\nu^4}
e_ke_me_j
\Biggr]
v^i,
\ee
where $\eta^i{}_j=a^{ih}\eta_{hj}$.
Noting that
\be
b_iv^i=(1-c^2)b, \qquad b_i\eta^i{}_j=-(1-c^2)e_j,
\ee
and
\be
u_i\eta^i{}_j=0, \qquad u_iv^i=q^2,
\ee
and taking into account the representation
$ y_i=(u_i+gqb_i) K^2 / B   $
(which is valid in the  space
 $\cF\cR^{PD}_{g;c}$
 under study
(see (A.27)),
we can readily conclude that
\be
y_iG^i{}_{kmj}=
-
(1-c^2)k\fr{g^2q^2}{\nu^2}
(\eta_{kj}e_m+\eta_{mj}e_k+\eta_{mk}e_j)  \fr{K^2}B.
\ee

{%\pgbrk}

Indeed, we have
$$
(u_i+gqb_i)G^i{}_{kmj}=
$$
\ses
\ses
$$
gq^2
\Biggl[
2k(1-c^2)^2\fr{g^2b}{\nu^3}
(\eta_{kj}e_m+\eta_{mj}e_k+\eta_{mk}e_j)
+6k(1-c^2)^3\fr{g^3q^2}{\nu^4}
e_ke_me_j
\Biggr]
$$

\ses
\ses
$$
-(1-c^2)g^2q(U_{kj} e_m   +   U_{mj} e_k   +U_{km} e_j)
$$
\ses
$$
+g^2qb(1-c^2)
\Biggl[
2k(1-c^2)^2\fr{g^2b}{\nu^3}
(\eta_{kj}e_m+\eta_{mj}e_k+\eta_{mk}e_j)
+6k(1-c^2)^3\fr{g^3q^2}{\nu^4}
e_ke_me_j
\Biggr]
$$

\ses

\ses

\ses

$$
=
gqk
\Biggl[
2(1-c^2)^2\fr{g^2b}{\nu^2}
(\eta_{kj}e_m+\eta_{mj}e_k+\eta_{mk}e_j)
+6(1-c^2)^3\fr{g^3q^2}{\nu^3}
e_ke_me_j
\Biggr]
$$
\ses
\ses
$$
-g^2kq(1-c^2)
\Biggl[
\Bigl(
\fr{2\nu\!-\!q}{\nu^2}\eta_{kj}
+2(1-c^2)^2\fr{g^2q^2}{\nu^3}e_ke_j
\Bigr)
e_m
+\Bigl(
\fr{2\nu\!-\!q}{\nu^2}\eta_{mj}
+2(1-c^2)^2\fr{g^2q^2}{\nu^3}e_me_j
\Bigr)
e_k
$$
\ses
$$
+\Bigl(
\fr{2\nu\!-\!q}{\nu^2}\eta_{km}
+2(1-c^2)^2\fr{g^2q^2}{\nu^3}e_ke_m
\Bigr)
e_j
\Biggr]
$$

\ses

\ses

\ses

$$
=
2gkq
(1-c^2)(\nu-q)\fr{g}{\nu^2}
(\eta_{kj}e_m+\eta_{mj}e_k+\eta_{mk}e_j)
$$
\ses
\ses
$$
-g^2kq(1-c^2)
\Biggl[
\fr{2\nu-q}{\nu^2}\eta_{kj}
e_m
+
\fr{2\nu-q}{\nu^2}\eta_{mj}
e_k
+
\fr{2\nu-q}{\nu^2}\eta_{km}
e_j
\Biggr],
$$
\ses
which shows that (B.14) is valid.

{%\pgbrk}

Let us also verify that (B.3) entails (B.10):
$$
U_{kmj}
+\fr1{q^2}(v_kU_{mj}+v_mU_{kj}+v_jU_{km})
$$

\ses

$$
=\fr2{q^2}
k\fr{\nu-q}{\nu^2}(v_j\eta_{km}  +  v_k\eta_{jm}  +  v_m\eta_{kj}  )
$$

\ses

$$
+2
k(1-c^2)^2\fr{g^2}{\nu^3}(e_ke_mv_j+  e_ke_jv_m+  e_je_mv_k)
$$

    \ses             \ses

$$
-6k\fr{q^2}{\nu^4}
\fr1{b^3}\Biggl[
q^3e_ke_me_j+q^2{\nu }(e_ke_mb_j+  e_ke_jb_m+e_je_mb_k)
+q{\nu }^2(b_kb_me_j+  b_kb_je_m+b_jb_me_k)
+{\nu }^3b_kb_mb_j
\Biggr]
$$

           \ses

$$
+4k\fr{q^2}{b^3\nu^3}
\Biggl[
[q^2e_ke_m+q{\nu }(e_kb_m+e_mb_k)+{\nu }^2b_kb_m]   e_j
+
[q^2e_ke_m+q{\nu }(e_kb_m+e_mb_k)+{\nu }^2b_kb_m]   b_j
\Biggr]$$

\ses

\ses

$$
+4k\fr{q^2}{b^3\nu^3}
\Biggl[
[q^2e_ke_j+q{\nu }(e_kb_j+e_jb_k)+{\nu }^2b_kb_j]   e_m
+
[q^2e_ke_j+q{\nu }(e_kb_j+e_jb_k)+{\nu }^2b_kb_j]   b_m
\Biggr]
$$

\ses
\ses

$$
+4k\fr{q^2}{b^3\nu^3}
\Biggl[
[q^2e_je_m+q{\nu }(e_jb_m+e_mb_j)+{\nu }^2b_jb_m]   e_k
+
[q^2e_je_m+q{\nu }(e_jb_m+e_mb_j)+{\nu }^2b_jb_m]   b_k
\Biggr]$$

               \ses
$$
+2k\fr{ q-\nu}{\nu^3}
(\nu_m\eta_{kj}   +   \nu_k\eta_{mj}    +   \nu_j\eta_{km})
$$

\ses

$$
-
2k\fr {q^2}{b^3\nu^2}
\Biggl[
q[e_ke_m+e_kb_m+e_mb_k+b_kb_m]   e_j
+
\nu [e_ke_m+e_kb_m+e_mb_k+b_kb_m]   b_j
\Biggr]
$$

    \ses

    \ses

$$
-
2k\fr {q^2}{b^3\nu^2}
\Biggl[
q[e_ke_j+e_kb_j+e_jb_k+b_kb_j]   e_m
+
\nu [e_ke_j+e_kb_j+e_jb_k+b_kb_j]   b_m
\Biggr]
$$

\ses

\be
-
2k\fr {q^2}{b^3\nu^2}
\Biggl[
q[e_je_m+e_jb_m+e_mb_j+b_jb_m]   e_k
+
\nu [e_je_m+e_jb_m+e_mb_j+b_jb_m]   b_k
\Biggr],
\ee
therefore  (B.10) is issued.

{%\pgbrk}

The representation (B.14) can be written as
\be
y_iG^i{}_{kmj}=
-
(1-c^2)k\fr{g^2q^2}{\nu^2}
(h_{kj}e_m+h_{mj}e_k+h_{mk}e_j)\
+
(1-c^2)k\fr{3g^2q^4}{B\nu^2}
e_ke_me_j
\fr{K^2}B,
\ee
\ses
where
we have used the following representation of the angular metric tensor:
\ses
\be
h_{ij}=
\fr{K^2}B
\lf(\eta_{ij}
+\fr{q^2}{B}e_ie_j
\rg)
\ee
(see (A.51) in [7]).

We can use  the equality
\be
A_i=-\fr {Kgq}{2B}
\fr1X
e_i
\ee
(compare (A.31) with (B.6))
and
 take into account the representation (A.34) for the Cartan tensor, which turns   (B.16)  to
\be
y_iG^i{}_{kmj}=
-
4(1-c^2)k
\lf(m_1A_{kmj}
+m_2A_kA_mA_j
\rg),
\ee
where
\be
m_1=      -        \fr{gqB}{2K\nu^2},
\qquad
m_2=  -  m_1 X\lf(N+1-\fr1X\rg)
\fr1{A_hA^h}.
\ee

The representation (B.19) tells us that
\ses
\be
\dot A_{kmj}=
(1-c^2)k
\bigl(m_1A_{kmj}
+m_2A_kA_mA_j
\bigr).
\ee
We have applied the known formula
$$\dot A_{kmj}
=-\fr14
y_iG^i{}_{kmj}
$$
(see p. 67 in [2]).

{%\pgbrk}

\vskip 1cm

\def\bibit[#1]#2\par{\rm\noindent\parskip1pt
                     \parbox[t]{.05\textwidth}{\mbox{}\hfill[#1]}\hfill
                     \parbox[t]{.925\textwidth}{\baselineskip11pt#2}\par}

\nin
{  REFERENCES}

\ses

\bibit[1] H. Rund: \it The Differential Geometry of Finsler
 Spaces, \rm Springer, Berlin 1959.

\ses

\bibit[2] D. Bao, S. S. Chern, and Z. Shen: {\it  An
Introduction to Riemann-Finsler Geometry,}  Springer, N.Y., Berlin 2000.

          %%%%%%%

\ses

\bibit[3] G. S. Asanov:  Finsleroid--Finsler  space with Berwald and  Landsberg conditions,
 {\it  arXiv:math.DG}/0603472 (2006).

\ses

\bibit[4] G. S. Asanov:  Finsleroid--Finsler  space and spray   coefficients, {\it  arXiv:math.DG}/0604526 (2006).

\ses

 \bibit[5] G. S. Asanov:  Finsleroid--Finsler  spaces of positive--definite and  relativistic types.
 \it Rep. Math. Phys. \bf 58 \rm(2006), 275--300.

\ses

\bibit[6] G. S. Asanov:  Finsleroid--Finsler space and geodesic spray    coefficients,
{\it Publ.  Math. Debrecen } {\bf 71/3-4} (2007), 397-412.

\ses

\bibit[7] G. S. Asanov:   Finsleroid-regular    space  developed.    Berwald case, {\it  arXiv:math.DG}/0711.4180v1 (2007).

\end{document}